\documentstyle[11pt,amssymb,amsfonts]{article}

\begin{document}
\newcommand{\p}{\parallel }
\makeatletter \makeatother
\newtheorem{th}{Theorem}[section]
\newtheorem{lem}{Lemma}[section]
\newtheorem{de}{Definition}[section]
\newtheorem{rem}{Remark}[section]
\newtheorem{cor}{Corollary}[section]
\renewcommand{\theequation}{\thesection.\arabic {equation}}

\title{{\bf A note on equivariant eta forms}
}
\author{Yong Wang \\}

\date{}
\maketitle

\begin{abstract} In this note, we prove the regularity of eta forms by the Clifford asymptotics. Then we generalize this
result to the equivariant case.\\

\noindent{\bf Keywords:}\quad
Eta forms; Equivariant eta forms; Clifford asymptotics\\
\end{abstract}

\section{Introduction}
    \quad In [APS], Atiyah, Patodi and Singer proved the Atiyah-Patodi-Singer index theorem
    for manifolds with boundary and they introduced the eta
    invariants. Bismut and Freed gave a simple proof of the regularity of
    eta invariants in [BF]. Bismut and Cheeger generalized the Atiyah-Patodi-Singer index
    theorem to the family case in [BC1,2]. They used the eta forms
    for families of Dirac operators. The regularity of eta forms was
    proved by the probability method in [BGS]. Donnelly generalized the Atiyah-Patodi-Singer index theorem
    to the equivariant case and introduced the equivariant eta
    invariants in [D]. Zhang proved the regularity of the equivariant eta
    invariants by the Clifford asymptotics in [Z2]. In this note, we
    firstly prove the regularity of eta forms by the Clifford
    asymptotics. Then we define the equivariant eta forms and prove
    their regularity.

\section{The regularity of eta forms
 }

 \quad Let $M$ be a $n+m$ dimensional compact connected manifold
 and $n$ be an odd integer and $B$ be a $m$ dimensional compact connected
 manifold. We assume that $\pi :M\rightarrow B$ is a submersion of
 $M$ onto $B$, which defines a fibration of $M$ with fibre $G$. For
 $y\in B$, $\pi^{-1}(y)$ is then a submanifolds $G_y$ of $M$. $TG$
 denotes the $n$-dimensional vector bundle on $M$ whose fibre $T_xG$
 is the tangent space at $x$ to the fibre $G_{\pi x}$. We assume
 that $M$ and $B$ are oriented. Taking the orthogonal bundle of $TG$
 in $TM$ with respect to any Riemannian metric, determines a smooth
 horizontal subbundle $T^HM$, i.e. $TM=T^HM\oplus TG$. Vector fields
 $X\in TB$ will be identified with their horizontal lifts $X\in
 T^HM$, moreover $T^H_xM$ is isomorphic to $T_{\pi(x)}B$ via
 $\pi_*$. Recall that $B$ is Riemannian, so we can lift the
 Euclidean scalar product $g_B$ of $TB$ to $T^HM$.
And we assume that $TG$ is endowed with a scalar product $g_G$. Thus
we can introduce in $TM$ a new scalar product $g_B\oplus g_G$, and
denote by $\nabla^L$ the Levi-Civita connection on $TM$ with respect
to this metric. Let $\nabla^B$ denote the Levi-Civita connection on
$TB$ and we still denote by $\nabla^B$ the pullback connection on
$T^HM$. Let $\nabla^G=P_G(\nabla^L)$ where $P_G$ denotes the
projection to $TG$. Let $\nabla^{\oplus}=\nabla^B\oplus \nabla^G$
and $S=\nabla^L-\nabla^{\oplus}$ and $T$ is the torsion tensor of
$\nabla^{\oplus}$. Let $SO(TG)$ be the $SO(n)$ bundle of oriented
orthonormal frames in $TG$. Now we assume that bundle $TG$ is spin.
Let $S(TG)$ be the associated spinors bundle and $\nabla^G$ can be
lifted to give a connection on $S(TG)$. Let $D$ be the tangent Dirac
operator. Let $K$ be the scalar curvature of fiber $G$ and
$e_1(x),\cdots e_n(x)$ denote the orthonormal frame of $TG$.
   If $A(Y)$ is any $0$ order operator depending linearly on
   $Y\in TM$, we define the operator $(\nabla_{e_i}+A(e_i))^2$ as
   follows
   $$(\nabla_{e_i}+A(e_i))^2=\sum_1^n(\nabla_{e_i(x)}+A(e_i(x)))^2-\nabla_{\sum_j\nabla_{e_j}e_j}
   -A(\sum_j\nabla_{e_j}e_j).\eqno(2.1)$$
Let ${\rm Tr}^{\rm even}$ denote taking trace on the coefficients of
even forms on $B$. Let $c(T)=\sum_{1\leq \alpha<\beta \leq
m}dy_\alpha dy_\beta c(T(\frac{\partial}{\partial
y_\alpha},\frac{\partial}{\partial y_\beta})),$ and
$$I^t=-t(\nabla^G_{e_i}+\frac{1}{2\sqrt{t}}<S(e_i)e_j,
f_\alpha>e_jdy_\alpha+\frac{1}{4t}<S(e_i)f_\alpha,f_\beta>dy_\alpha
dy_\beta)^2+\frac{tK}{4}.\eqno(2.2)$$ We make the following
definition of the eta forms.\\

\noindent {\bf Definition 2.1} $\hat{\eta}$ denotes the even degree
form on $B$
$$\hat{\eta}=\frac{1}{2\sqrt{\pi}}\int _0^{+\infty}{\rm Tr}^{\rm
even}[(D+\frac{c(T)}{4t}){\rm
exp}(I^t)]dt/t^{\frac{1}{2}}.\eqno(2.3)$$\\

\indent In order to prove that (2.3) is well defined, we need to
prove the regularity at origin and $+\infty$. Firstly we prove the
regularity at $+\infty$. We consider the two terms ${\rm Tr}^{\rm
even}[D{\rm exp}(I^t)] $ and ${\rm Tr}^{\rm
even}[\frac{c(T)}{4t}{\rm exp}(I^t)]$. For the first term,
considering $D$, we can assume that $D$ is invertable and using the
discussions as in [W, pp.148-150], then we can get the regularity at
$+\infty$ of the first term. Since $\frac{c(T)}{4t}$ is a bounded
operator, using the discussions in [BGV, p. 275], we can get the
regularity at $+\infty$ of the second term. Nextly, we prove the
regularity at origin. Fixing a $y\in B$ and considering the fibre
$G_y$, let ${\rm exp}(I^t)(x,x')$ for $x,x'\in G_y$ denote the heat
kernel of $I^t$ (see [BGV]). By a simple discussion, we can get
$${\rm Tr}^{\rm even}[(D+\frac{c(T)}{4t}){\rm exp}(I^t)]=\int_{G_y}
{\rm Tr}^{\rm even}[(D+\frac{c(T)}{4t}){\rm
exp}(I^t)(x,x)]dx.\eqno(2.4)$$ Let $I^1=I^t|_{t=1}$. By the Duhamel
principle, we can get
$${\rm Tr}^{\rm even}[{\rm exp}(t(I^1+dt(D+\frac{c(T)}{4})))]=
{\rm Tr}^{\rm even}[{\rm exp}(tI^1)]+tdt{\rm Tr}^{\rm
even}[(D+\frac{c(T)}{4}){\rm exp}(tI^1)].\eqno(2.5)$$
\quad Let
$\phi$ be a smooth function defined locally in a neighborhood of
$x$, denote the degree of zero of $\phi$ at $x$ by $\upsilon(\phi)$,
to every
$$\alpha(x')=\phi_{l_1}(x')\frac{\partial}{\partial x_{i_1}}
\phi_{l_2}(x')\cdots\frac{\partial}{\partial x_{i_m}}
\phi_{l_{m+1}}(x')dy_{\alpha_1}\cdots dy_{\alpha_p}dt$$
$$e_{j_1}\cdots e_{j_s}: S(TG)_{x'}\rightarrow S(TG)_{x'};
\alpha_i\neq \alpha_j~ (i\neq j);~ j_a\neq j_t~ (a\neq
t),\eqno(2.6a)$$ we define
$$
\chi(\alpha)=m+p+s+1-\upsilon(\phi_1\cdots\phi_{m+1}).\eqno(2.6b)$$
We call $\alpha(x')$ an even (odd) element if $p+s$ is a even (odd)
integer and we denote $\{\chi<m\}$ the linear space generated by all
the elements $\alpha$ for which $\chi(\alpha)<m$ and denote
$(\chi<m)$ an element of $\{\chi<m\}$, e.g.
$\omega=\omega'+(\chi<m)$ means that there exists a $\beta\in
\{\chi<m\}$ such that $\omega=\omega'+\beta$. Set
$$h(x)=1+\frac{1}{2}dt\sum_{i=1}^nx_ie_i.\eqno(2.7)$$
Then we have
$$ he_ih^{-1}=e_i+dt(\chi\leq -1);~~
h(\frac{1}{2}e_idt)h^{-1}=\frac{1}{2}e_idt,\eqno(2.8)$$
$$~h\nabla^G_{e_i}h^{-1}=\nabla^G_{e_i}-\frac{1}{2}dte_i+dt(\chi\leq
-1);~ \eqno(2.9)$$
$$h(\frac{1}{2}<S(e_i)e_j,
f_\alpha>e_jdy_\alpha )h^{-1}=\frac{1}{2}<S(e_i)e_j,
f_\alpha>e_jdy_\alpha +dt(\chi\leq -1);~\eqno(2.10)$$
$$h(\frac{1}{4}<S(e_i)f_\alpha,f_\beta>dy_\alpha
dy_\beta)h^{-1}=\frac{1}{4}<S(e_i)f_\alpha,f_\beta>dy_\alpha
dy_\beta.\eqno(2.11)$$ By the proposition 2.10 in [BGS], we have
$$I^1+dt(D+\frac{c(T)}{4})=-(\nabla^G_{e_i}+\frac{1}{2}<S(e_i)e_j,
f_\alpha>e_jdy_\alpha$$
$$+\frac{1}{4}<S(e_i)f_\alpha,f_\beta>dy_\alpha
dy_\beta-\frac{1}{2}e_idt)^2+\frac{K}{4},\eqno(2.12)$$ By
(2.7)-(2.12), we have
$$h[I^1+dt(D+\frac{c(T)}{4})]h^{-1}=I^1+dtu,\eqno(2.13)$$
where $\chi(u)\leq 0$. For $t>0$, by [Z1,(4.17)], we have
$${\rm exp}(tI^1)(x,x')=\frac {e^{-d(x,x')^2/4t}}{(4\pi
t)^{\frac{n}{2}}}\left(\sum_{i=0}^{[\frac{n}{2}]+[\frac{m}{2}]+2}U_it^i+o(t^{[\frac{n}{2}]
+[\frac{m}{2}]+2})\right). \eqno(2.14)$$ For $I^1+dtu$, similar to
(4.17) in [Z1], we have
$$[\hat{d}+x_i(b_i\Gamma^\alpha_{ij}e_jdy_\alpha+c_i\Gamma^\beta_{i\alpha}dy_\alpha
dy_\beta+B_i)+r]\widetilde{U_r}+dtA\widetilde{U_r}=(I^1+dtu)\widetilde{U_{r-1}}\eqno(2.15)$$
where $\widetilde{U_r}=U_r+dtV_r$ and $\widetilde{U_{-1}}=0$ and
$\chi(A)\leq -2$ and see [Z1] for
$\hat{d},~b_i\Gamma^\alpha_{ij},~B_i$. By (2.15), we have
$${\rm exp}(t(I^1+dtu))(x,x')=\frac {e^{-d(x,x')^2/4t}}{(4\pi
t)^{\frac{n}{2}}}\left(\sum_{i=0}^{[\frac{n}{2}]+[\frac{m}{2}]+2}(U_i+dtV_i)t^i+o(t^{[\frac{n}{2}]
+[\frac{m}{2}]+2})\right), \eqno(2.16)$$ where $\chi(U_i)\leq 2i$,
$\chi(V_i)\leq 2(i-1)$ and $U_i, V_i$ contain no $dt$. $U_i$ is an
even element and $V_i$ is an odd element. By (2.13), we have
$${\rm exp}[t(I^1+dt(D+\frac{c(T)}{4}))](x,x)=h^{-1}(x){\rm exp}[t(I^1+dtu)]h(x).\eqno(2.17)$$
By (2.5) and (2.17), we get
$$tdt{\rm Tr}^{\rm
even}[(D+\frac{c(T)}{4}){\rm exp}(tI^1)(x,x)]={\rm Tr}^{\rm
even}\{h^{-1}(x){\rm exp}[t(I^1+dtu)(x,x)]h(x)\}$$ $$-{\rm Tr}^{\rm
even}[{\rm exp}(tI^1)(x,x)].\eqno(2.18)$$ By (2.8), (2.14) and
(2.16), we get $${\rm Tr}^{\rm even}\{h^{-1}(x){\rm
exp}[t(I^1+dtu)(x,x)]h(x)\}-{\rm Tr}^{\rm even}[{\rm
exp}(tI^1)(x,x)]$$ $$=\frac {1}{(4\pi t)^{\frac{n}{2}}}{\rm Tr}^{\rm
even}\left(\sum_{i=0}^{[\frac{n}{2}]+[\frac{m}{2}]+2}dtW_it^i+o(t^{[\frac{n}{2}]
+[\frac{m}{2}]+2})\right), \eqno(2.19)$$ where $\chi(W_i)\leq
2(i-1)$ is an odd element. Let $\phi_t$ be the rescaling operator
defined by
$\phi_t(dy_\alpha)=\frac{1}{\sqrt{t}}dy_\alpha,~\phi_t(dt)=\frac{1}{\sqrt{t}}
dt.$ Then by (2.18) and (2.19), we have $${\rm Tr}^{\rm
even}[(D+\frac{c(T)}{4t}){\rm exp}(I^t)(x,x)]=\frac {1}{(4\pi
t)^{\frac{n}{2}}t}\phi_t\left\{{\rm Tr}^{\rm
even}\left(\sum_{i=0}^{[\frac{n}{2}]+[\frac{m}{2}]+2}W_it^i+o(t^{[\frac{n}{2}]
+[\frac{m}{2}]+2})\right)\right\}. \eqno(2.20)$$\\

\noindent {\bf Lemma 2.2}~ {\it Suppose $i\leq [\frac{n}{2}]+1.$ If
$W$ is an odd element and $\chi(W_i)\leq 2(i-1)$, then ${\rm
Tr}^{\rm
even}W_i(0)=0$}\\
\noindent{\it Proof.} Since we take the even form trace, we can
assume $$W_i=cx_{i_1}\cdots x_{i_k}e_{j_1}\cdots
e_{j_l}dy_{\beta_1}\cdots dy_{\beta_{2s}}.$$ If $k>0$, we have
$W_i(0)=0$. When $k=0$, $\chi(W_i)\leq 2(i-1)\leq
2([\frac{n}{2}]+1-1)=n-1.$ By $\chi(dy_\beta)>0$ and $W_i$ being an
odd element, we have $l<n$ is an odd integer, so tr$(e_{j_1}\cdots
e_{j_l})=0$.~~~~$\Box$\\

\noindent {\bf Lemma 2.3}~ {\it Suppose $1\leq j\leq
[\frac{m}{2}]+1.$ If $W_{[\frac{n}{2}]+1+j}$ is an odd element and
$\chi(W_{[\frac{n}{2}]+1+j})\leq 2([\frac{n}{2}]+j)$, then}
$$\phi_t\left(\frac{1}{t^{\frac{n}{2}+1}}{\rm Tr}^{\rm
even}W_{[\frac{n}{2}]+1+j}(0)t^{[\frac{n}{2}]+1+j}\right)=O(t^{\frac{1}{2}}),~~~~~~~t\searrow 0.\eqno(2.21)$$\\
\noindent{\it Proof.}We can assume
$$W_{[\frac{n}{2}]+1+j}=cx_{i_1}\cdots x_{i_k}e_{j_1}\cdots
e_{j_l}dy_{\beta_1}\cdots dy_{\beta_{2s}}.$$ Then $k=0$ and $l=n$,
otherwise ${\rm Tr}^{\rm even}W_{[\frac{n}{2}]+1+j}(0)=0.$ Since
$\chi(W_{[\frac{n}{2}]+1+j})\leq 2([\frac{n}{2}]+j)=n-1+2j$, then
$s\leq j-1$ and
$\phi_t(W_{[\frac{n}{2}]+1+j})=t^{-s}W_{[\frac{n}{2}]+1+j}$. So the
degree of $t$ is $[\frac{n}{2}]+1+j-(\frac{n}{2}+s+1)\geq
\frac{1}{2}$.~~~~~~~~$\Box$\\

By (2.20), Lemmas 2.2 and 2.3, we get\\

\noindent {\bf Theorem 2.4([BGS])}~
$${\rm Tr}^{\rm
even}[(D+\frac{c(T)}{4t}){\rm
exp}(I^t)(x,x)]=O(t^{\frac{1}{2}}),~~~~~~~t\searrow 0.\eqno(2.22)$$

\section { The regularity of equivariant eta forms}

  \quad Let isometry $g$ act fibrewise on $M$ and act as identity
  on $B$ and $g$ preserve the orientation and the spin structure on
  $S(TG)$. Let $\widetilde{g}:\Gamma (S(TG))\rightarrow \Gamma (S(TG))$
  be the lift of $g$. We have the following definition of equivariant eta forms.\\

  \noindent{\bf Definition 3.1}
$\hat{\eta}_g$ denotes the even degree form on $B$
$$\hat{\eta}_g=\frac{1}{2\sqrt{\pi}}\int _0^{+\infty}{\rm Tr}^{\rm
even}[\widetilde{g}(D+\frac{c(T)}{4t}){\rm
exp}(I^t)]dt/t^{\frac{1}{2}}.\eqno(3.1)$$\\

\indent Since $\widetilde{g}$ is bounded operator, we can prove the
regularity at $+\infty$ as in Section 2. Nextly, we prove the
regularity at origin. Fix a $y\in B$ and consider the fibre $G_y$.
Similar to Corollary 1.4 in [Z2], if $g$ has no fixed points on
$G_y$, then we can prove
$$||\widetilde{g}(D+\frac{c(T)}{4t}){\rm exp}(I^t)(x,gx)||\leq
\frac{C_1}{t^{\frac{n}{2}+\frac{m}{2}+1}}{\rm
exp}(-C_2/t),\eqno(3.2)$$ where $C_1, C2$ are constants and the norm
represents the norm of coefficients of forms on $B$. So in this
case, (3.1) is well defined. Since $g$ is an isometry on $G_y$, for
convenience, we assume the fixed point set $F$ of $g$ is connected
and ${\rm codim}F=2n'$ and ${\rm dim}F=q$. Denote by $N(F)$ the
normal bundle to $F$. Similar to the discussions of (3.2), we need
only to prove that
$$ {\rm lim}_{t\rightarrow 0}\left|\int_F\int_{N_\xi(\varepsilon)}{\rm Tr}^{\rm
even}[\widetilde{g}(D+\frac{c(T)}{4t}){\rm
exp}(I^t)(x,gx)]dN_{\xi}d\xi\right|\leq C,\eqno(3.3)$$ for some
constant $C>0$. Here $N_\xi(\varepsilon)=\left\{v\in
N_\xi(F)|||v||<\varepsilon\right\}$. Similar to (2.18), we have
$$tdt{\rm Tr}^{\rm
even}[\widetilde{g}(D+\frac{c(T)}{4}){\rm exp}(tI^1)(x,gx)]={\rm
Tr}^{\rm even}\{\widetilde{g}h^{-1}(x){\rm
exp}[t(I^1+dtu)(x,gx)]h(gx)\}$$
$$-{\rm Tr}^{\rm even}[\widetilde{g}{\rm exp}(tI^1)(x,gx)].\eqno(3.4)$$
By (2.14),(2.16) and (3.4), we have
$$t\widetilde{g}(D+\frac{c(T)}{4}){\rm exp}(tI^1)(x,gx)~~~~~~~~~~~~~~~~~~~$$
$$= \frac {e^{-d(x,gx)^2/4t}}{(4\pi
t)^{\frac{n}{2}}}\widetilde{g}\left\{
\sum_{i=1}^{2n'}(((dg-I)x)_{q+i}e_{q+i})\left(\sum_{j=0}^{[\frac{n}{2}]+[\frac{m}{2}]+2}
U_jt^j+o(t^{[\frac{n}{2}] +[\frac{m}{2}]+2})\right)\right.$$
$$~~~~~~~~~~~~~~~~~~~~~~~~~\left.
+\sum_{i=0}^{[\frac{n}{2}]+[\frac{m}{2}]+2}W_it^i+o(t^{[\frac{n}{2}]
+[\frac{m}{2}]+2})\right\}, \eqno(3.5)$$ where $\chi(W_i)\leq
2(i-1)$ is an odd element. By [Z2, p.1126], we know that
$$\chi(\widetilde{g})\leq 2n',~~~\chi(\widetilde{g}\sum_{i=1}^{2n'}(((dg-I)x)_{q+i}e_{q+i}))\leq
2n'-2. \eqno(3.6)$$

\noindent {\bf Lemma 3.2}~ {\it Suppose $1\leq j\leq
[\frac{n}{2}]+[\frac{m}{2}]+2.$ If $\overline{W}$ is an odd element
and $\chi(\overline{W})\leq 2n'+2j-2$,  then}
$${\rm lim}_{t\rightarrow
0}(\frac{1}{t})^{3/2}\left|\int_{N_\xi(\varepsilon)}\frac
{e^{-d(x,gx)^2/4t}}{(4\pi t)^{\frac{n}{2}}}\phi_t[{\rm
tr}(\overline{W}(0;x)]t^jdx\right|\leq C_1,\eqno(3.7)$$ {\it for
some constant $C_1>0$; where in the $W(z;x)$. $z$ stands
for tangential coordinates and $x$ stands for normal coordinates.}\\
\noindent{\it Proof.} We can assume that $\overline{W}$ is a
monomial, then it can be written as
$$\overline{W}=\phi(0)x_{i_1}\cdots x_{i_k}e_1\cdots
e_ndy_{\beta_1}\cdots dy_{\beta_{2s}}.$$ We note that we can assume
that $x_i$ in $\overline{W}$ are normal coordinates, for otherwise
${\rm tr}\overline{W}(0;\dot)=0.$\\
\indent (i) If $\chi(\overline{W})=2n'+2j-2,$ then $k=n+2s-2n'-2j+2$
is an odd integer. By making the change of variables $x=t^{1/2}b$,
we get
$${\rm lim}_{t\rightarrow
0}(\frac{1}{t})^{3/2}\left|\int_{N_\xi(\varepsilon)}\frac
{e^{-d(x,gx)^2/4t}}{(4\pi t)^{\frac{n}{2}}}\phi_t[{\rm
tr}(\overline{W}(0;x)]t^jdx\right|$$
$$\leq {\rm lim}_{t\rightarrow
0}(\frac{1}{t})^{3/2}\left|\int_{N_\xi(\varepsilon/\sqrt{t})}\frac
{e^{-||(I-dg)b||^2}}{(4\pi
t)^{\frac{n}{2}}}t^{\frac{n}{2}+s-n'-j+1}b_{i_1}\cdots
b_{i_k}t^j\frac{1}{t^s}dy_{\beta_1}\cdots
dy_{\beta_{2s}}t^{n'}db\right|=0$$ \indent (ii) If
$\chi(\overline{W})<2n'+2j-2,$ then $k>n+2s-2n'-2j+2$ and we get
$${\rm lim}_{t\rightarrow
0}(\frac{1}{t})^{3/2}\left|\int_{N_\xi(\varepsilon)}\frac
{e^{-d(x,gx)^2/4t}}{(4\pi t)^{\frac{n}{2}}}\phi_t[{\rm
tr}(\overline{W}(0;x)]t^jdx\right|\leq C_3.$$ We note that the above
discussions are also correct when $n+2s-2n'-2j+2<0$.~~$\Box$
\indent By (3.5), (3.6) and Lemma 3.2, we get\\

 \noindent {\bf Theorem 3.3}~
$${\rm Tr}^{\rm
even}[\widetilde{g}(D+\frac{c(T)}{4t}){\rm
exp}(I^t)(x,gx)]=O(t^{\frac{1}{2}}),~~~~~~~t\searrow
0.\eqno(3.8)$$\\

 \noindent {\bf Acknowledgement.} This work
was supported by NSFC No.10801027 and Fok Ying Tong Education
Foundation No. 121003.\\

\noindent{\large \bf References}\\

\noindent[APS]Atiyah, M. F.; Patodi, V. K.; Singer, I. M., Spectral
asymmetry and Riemannian geometry. I. Math. Proc. Cambridge Philos.
Soc. 77(1975), 43-69. \\
\noindent[BGV]Berline, N.; Getzler, E.; Vergne, M.,{\it Heat kernels
and Dirac operators.} Springer-Verlag, Berlin, 1992.\\
\noindent[BC1]Bismut, J. M.; Cheeger, J., Families index for
manifolds with boundary, superconnections, and cones. I. Families of
manifolds with boundary and Dirac operators. J. Funct. Anal. 89
(1990), no. 2,
313-363.\\
\noindent[BC2]Bismut, J. M.; Cheeger, J., Families index for
manifolds with boundary, superconnections and cones. II. The Chern
character. J. Funct. Anal. 90 (1990), no. 2, 306-354.\\
\noindent[BF]Bismut, J. M.; Freed, D., The analysis of elliptic
families. II. Dirac operators, eta invariants, and the holonomy
theorem. Comm. Math. Phys. 107 (1986), no. 1, 103-163.\\
\noindent[BGS]Bismut, J.-M.; Gillet, H.; Soul¨¦, C. Analytic torsion
and holomorphic determinant bundles. I. Bott-Chern forms and
analytic torsion. Comm. Math. Phys. 115 (1988), no. 1, 49-78.\\
\noindent[D]Donnelly, H., Eta invariants for $G$-spaces. Indiana
Univ. Math. J. 27 (1978), no. 6, 889-918.\\
\noindent[W]Wu, F., The Chern-Connes character for the Dirac
operator on manifolds with boundary. $K$-Theory 7 (1993), no. 2,
145-174.\\
 \noindent[Z1]Zhang, W. P., Local Atiyah-Singer index
theorem for families of Dirac operators. Differential geometry and
topology ,Lecture Notes in Math., 1369, Springer, Berlin, 1989 351-366.\\
\noindent[Z2]Zhang, W. P., A note on equivariant eta invariants.
Proc. Amer. Math. Soc. 108 (1990), no. 4, 1121-1129.\\

 \indent{  School of Mathematics and Statistics,
Northeast Normal University, Changchun Jilin, 130024, China }\\
\indent E-mail: {\it wangy581@nenu.edu.cn}\\

\end{document}